\documentclass[12pt,a4paper]{article}
\usepackage{amsfonts,amssymb,graphicx}
\usepackage{theorem}
\newcommand{\bbf}{{\bf b}}
\newcommand{\cbf}{{\bf c}}

\newcommand{\Xbf}{{\bf X}}
\newcommand{\Ybf}{{\bf Y}}

\newcommand{\Pcal}{{\cal P}}

\newcommand{\Tcal}{{\cal T}}
\newcommand{\Ucal}{{\cal U}}
\newcommand{\Vcal}{{\cal V}}

\newcommand{\N}{{\mathbb N}}
\newcommand{\chara}{{\rm char}\mbox{$\,$}}
\newcommand{\mod}{{\hspace{-1ex}}\pmod}
\newcommand{\ncal}[1]{\mbox{$ \cal N$$^{(#1)}$}}

\newtheorem{thm}{Theorem}
 {\theorembodyfont{\rm}
\newtheorem{exa}{Example}
\newtheorem{defi}{Definition}
\newtheorem{rem}{Remark}
 }
\title{
Veronese Varieties over Fields with non-zero Characteristic: A Survey
}
\author{Hans Havlicek}
\date{}

\begin{document}
\maketitle

\section{Introduction}\label{Kap-Intro}
Non-zero characteristic of the (commutative) ground field $F$ heavily
influences the geometric properties of Veronese varieties and, in
particular, normal rational curves. Best known is probably the fact
that, in case of characteristic two, all tangents of a conic are
concurrent. This has lead to the concept of a {\em nucleus}. However,
it seems that there are essentially distinct definitions. Some
authors, like J.A.\ Thas \cite{thas-69}, use this term to denote a
point which extends a normal rational curve to an $(q+2)$-arc ($F$ a
finite field of even order $q$), others, like A.\ Herzer
\cite{herz-82}, use the same term for the intersection of all
osculating hyperplanes of a Veronese variety. In order to overcome
this difference of terminology we introduce the term {\em
$(r,k)$-nucleus}. The two types of nuclei mentioned above are just
particular examples fitting into this general concept.

Each nucleus is an {\em invariant subspace}, i.e.\ a subspace in the
ambient space of a Veronese variety which is fixed (as a set of
points) under the group of automorphic collineations of the variety.
However, an invariant subspace needs not be a nucleus.

In the present survey we collect some recent results on nuclei of
Veronese varieties and invariant subspaces of normal rational curves.
We must assume, however, that the ground field is not ``too small'',
since otherwise a Veronese variety is like dust: ``few points'' in
some ``high-dimensional'' space.

Nuclei and invariant subspaces do not appear in classical
textbooks on Veronese varieties ($F={\mathbb R},{\mathbb C})$,
since for characteristic zero all invariant subspaces are trivial.
If the ground field has characteristic $p>0$, then geometric
properties of invariant subspaces are closely related to {\em
multinomial coefficients} that vanish modulo $p$\/ and to the
representations of certain integers in base $p$. In order to
illustrate this connection some results on binomial and
multinomial coefficients are gathered in Chapters \ref{Kap-Binom}
and \ref{Kap-Multinom}.

\section{Pascal's triangle modulo a prime $p$}\label{Kap-Binom}

\subsection{A partition of zero entries}

Throughout this section let $p$ be a fixed prime. The representation
of a non--negative integer $n\in\N:=\{0,1,2,\ldots\}$ in base $p$ has
the form
   \begin{equation}
   n = \sum\limits_{\sigma=0}^\infty n_\sigma p^\sigma
   =: \langle n_\sigma\rangle
   \end{equation}
with only finitely many digits $n_\sigma\in\{0,1,\ldots,p-1\}$
different from $0$.

Let $\langle n_\sigma\rangle$ and $\langle j_\sigma\rangle$ be the
representations of non--negative integers $n$ and $j$ in base $p$. By
a Theorem of Lucas \cite[364]{brou+w-95},
   \begin{equation}\label{lucas}
   {n \choose j} \equiv
   \prod_{\sigma=0}^\infty {n_\sigma \choose j_\sigma} \mod p.
   \end{equation}
{\em Pascal's triangle modulo} $p$\/ will be denoted by $\Delta$. The
numbering of its rows starts with the index $0$. Also, let $\Delta^i$
($i\in\N$) be the subtriangle of $\Delta$ that is formed by the rows
$0,1,\ldots, p^i-1$. From (\ref{lucas}) each triangle $\Delta^{i+1}$
($i\geq 0$) has the following form, with products taken modulo $p$\/:
   \renewcommand{\arraystretch}{1.}%
   \begin{displaymath}
   \begin{array}{c@{}c@{}c@{}c@{}c@{}c@{}c}
   {}& {}& {}&
   {0 \choose 0}\Delta^{i} &
   {}& {} & {}
   \\
   {}& {}&
   {1 \choose 0}\Delta^{i} &
   \nabla^{i} &
   {1 \choose 1}\Delta^{i}
   & {}& {}
   \\
   {}&
   {2 \choose 0}\Delta^{i}&
   \nabla^{i} &
   {2 \choose 1}\Delta^{i}&
   \nabla^{i}&
   {2 \choose 2}\Delta^{i} & {}
   \\
   \multicolumn{7}{c}{\quad\quad\dotfill\quad\quad}\\
   {p-1 \choose 0}\Delta^{i}\;\;\nabla^{i}&
   {} & {}& \ldots & {}& {}&
   \nabla^{i}\;\; {p-1 \choose p-1}\Delta^{i}
   \end{array}
   \end{displaymath}
   \renewcommand{\arraystretch}{1.0}%
Here the $\nabla^i$'s are triangles with all entries equal to zero.
Observe that the baseline of $\Delta^i$ has $p^i$ entries, whereas
the top line of $\nabla^i$ is formed by $p^i-1$ zero entries. So
$\nabla^0$ is empty. The binomial coefficients on the left hand side
of the $\Delta^i$'s are exactly the entries of $\Delta^1$.  None of
them vanishes modulo $p$. If $i\geq 2$, then each subtriangle
${n\choose j}\Delta^i$ from above can be decomposed into subtriangles
proportional to $\Delta^{i-1}$ and subtriangles $\nabla^{i-1}$, and
so on. Cf., among others, \cite[91--92]{hexe+s-78}, \cite[Theorem
1]{long-81b}, \cite{wolf-84}.

The zero entries of Pascal's triangle modulo $p$\/ fall into
(disjoint) maximal subtriangles $\nabla^i$ ($i\in \N^{+}$). We get a
partition of all zero entries of $\Delta$ by gluing together all
triangles $\nabla^i$ of same size to one class, say $\overline{i}$. A
formal definition of this partition, which is the backbone of many
further considerations, is as follows:
\begin{defi}\cite{gmai+h-00a}
A pair $(n,j)=(\langle n_\sigma\rangle,\langle j_\sigma\rangle)$ of
non--negative integers with $j\leq n$, ${n \choose j} \equiv 0 \mod
p$, and $L := \max \{ \sigma \in\N \mid j_\sigma > n_\sigma\}$, is in
{\em class}\/ $\overline{i}$, if
   \begin{equation}
   i = \min \{ \sigma \mid \sigma>L,\; j_\sigma < n_\sigma \}
   \in\N^{+}.
   \end{equation}
\end{defi}

\subsection{Counting zero entries}

The following is taken from \cite{gmai+h-00a}. Let $n=\langle
n_\sigma\rangle\in\N$ and $i\in\N^{+}$. Then the number of entries in
row $n$ of $\Delta$ belonging to class $\overline{i}$ equals
\begin{equation}\label{def-Phi(i,n)}
   \Phi(i,n):= \# \overline{i(n)}  =
   \Big( p^{i}-1- \sum_{\mu=0}^{i-1} n_\mu p^\mu \Big)\cdot n_{i} \cdot
   \prod_{\sigma=i+1}^\infty (n_\sigma +1).
\end{equation}
The number of entries in row $n$ of $\Delta$ belonging to classes
$\overline{i}$, $\overline{(i+1)}$, $\ldots$ is
   \begin{eqnarray}\label{Sigma(i,n)}
   \Sigma(i,n)  :=  \sum\limits_{\eta=i}^\infty \Phi(\eta,n)
    =  n + 1 - \Big(1 + \sum\limits_{\mu=0}^{i-1} n_\mu p^\mu\Big)
   \prod\limits_{\sigma=i}^\infty (n_\sigma + 1).
   \end{eqnarray}
For $i=1$ this is due to N.J.\ Fine \cite{fine-47}.

When exhibiting ``vertical'' properties of $\Delta$ the following
{\em top line function\/} turns out useful: Given $b\in\N^{+}$ and
$R\in\N$ then let
   \begin{equation}\label{T(i,b)}
   T(R,b):=  \sum_{\sigma=R}^{\infty} b_{\sigma}p^{\sigma}.
   \end{equation}
This function has the following property: If $(n,j)\in\overline i$
and $b:=n+1$ then $T(i,b)$ gives the ``top line'' of the triangle
$\nabla^i$ containing the $(n,j)$-entry of $\Delta$, i.e.\
   \begin{equation}
   0\equiv {n\choose j} \equiv {n-1\choose j} \equiv \ldots
   \equiv{T(i,b)\choose j}\not\equiv {T(i,b)-1\choose j} \mod p.
   \end{equation}

We refer to \cite{harb-75}, \cite{kara-96}, \cite{robe-57} for
further properties of $\Delta$.

\section{Normal rational curves}\label{Kap-NRC}

\subsection{Definition of $k$-nuclei}

Let $\{\bbf_0,\bbf_1\}$ be a basis of a $2$-dimensional vector space
$\Xbf$ over a commutative field $F$ (the parameter space) and let
$\Ybf$ be an $(n+1)$-dimensional vector space over $F$ with a basis
$\{\cbf_{0},\cbf_1,\ldots,\cbf_n\}$, where $n\geq 2$. The {\em
Veronese mapping}
   \begin{equation}\label{g-koo-NRC}
   F(x_0 \bbf_0+ x_1 \bbf_1)
   \mapsto
   F\Big( \sum_{e=0}^{n} x_0^{n-e} x_1^{e} \cbf_{e} \Big)
   \;\;(x_i\in F)
   \end{equation}
maps the point set of the projective line $\Pcal(\Xbf)$ into the
point set of $\Pcal(\Ybf)$, i.e.\ the projective space on $\Ybf$. Its
image is a {\em normal rational curve} $\Vcal_1^n$ (sometimes
abbreviated as NRC) with ambient space $\Pcal(\Ybf)$ \cite{bert-07},
\cite{bert-24}, \cite{brau-762}, \cite{bura-61}, \cite[Chapter
21]{hirs-85}. In terms of coordinates and an inhomogeneous parameter
$x:=x_1 / x_0$ we obtain
   \begin{equation}\label{param}
   \Vcal_1^n:=\{F(1,x,\ldots,x^n)\mid x\in F\cup\{\infty\} \}.
   \end{equation}

Recall the {\em non-iterative derivation} due to H.~Hasse,
F.K.~Schmidt, and O.~Teich\-m\"uller \cite{hass-37},
\cite[1.3]{hirs-98}. The $k$-th derivative
$D^{(k)}\,:\,F[X]\rightarrow F[X]$ is a linear mapping such that
$D^{(k)}(X^r)= {r\choose k}X^{r-k}$ for $k\leq r$ and $D^{(k)}(X^r) =
0$ otherwise ($r,k\in \N$).

If we fix one $u\in F$ then columns of the regular matrix
   \begin{equation}\label{C_u}
\left(
\renewcommand{\arraystretch}{1.4}
\begin{array}{ccccc}
{0\choose 0}   &0                  &0 &\ldots&0\\ {1\choose
0}u  &{1\choose 1}       &0 &\ldots&0\\ {2\choose
0}u^2&{2\choose 1}u      &{2\choose 2} &\ldots&0\\ \vdots &
&\ &\ddots&\vdots\\ {n\choose 0}u^n&{n\choose
1}u^{n-1}&{n\choose 2}u^{n-2}&\ldots&{n\choose n}
\end{array}
\right)
\renewcommand{\arraystretch}{1.0}
   \end{equation}
give, respectively, a point of the NRC (\ref{param}) and its {\em
derivative points}. The {\em $k$-osculating subspace} $(k\in\{-1,0,
\ldots, n-1\})$ of $\Vcal_1^n$ at the given point is the
$k$-dimensional projective subspace spanned by the first $k+1$
columns of the matrix (\ref{C_u}). The derivative points at $F\cbf_n$
$(u=\infty)$ are $F\cbf_{n-1},F\cbf_{n-2} \ldots,F\cbf_{0}$.

Formal derivation in $F[X]$ is in general not an adequate tool to
describe osculating subspaces \cite{ries-81}.

The (empty) $(-1)$-osculating subspace is introduced for formal
reasons only. However, we refrain from calling the entire space an
$n$-osculating subspace. As usual, a $1$-osculating subspace is also
called a {\em tangent}.

\begin{defi}\cite{gmai+h-00a}\label{def-k-knoten}
The $k$-{\em nucleus} $\ncal{k}\Vcal_1^n$ ($k\in\{-1,0,\ldots,n-1\}$)
of a normal rational curve $\Vcal_1^n$ is the intersection of all its
$k$-osculating subspaces.
\end{defi}

\begin{rem}\label{def-osk}
Instead of a parametric representation one could also use a {\em
generating map} \cite{havl-83}, \cite{havl-85}, {\em Segre varieties}
\cite{bura-61}, \cite{karz-87}, \cite{zeug-72}, \cite{zeug-77} or
tools from multilinear algebra \cite{gmai+h-00b}, \cite{herz-82} in
order to define osculating subspaces.

The NRC (\ref{param}) has exactly $\# F+1$ points. Every subset with
$s\leq n+1$ points is linearly independent. If $\# F\geq n+2$ then
the NRC is a set with at least $n+3$ points of which every $n+1$ are
linearly independent.

Suppose now that $q:=\# F$ is finite. From above, for $q \geq n+2$
the NRC is a non-trivial example of a $(q+1)$-arc in the
$n$-dimensional projective space $\Pcal(\Ybf)$. If $q=n+1$ then the
NRC is a frame. Furthermore, $u^n=u^{q-1}=1$ for all non-zero
elements $u$ of $F$. This illustrates that here the hyperplane with
equation $x_0=x_n$ contains all points of the NRC other than
$F\cbf_0$ and $F\cbf_n$. Finally, if $q \leq n$ then the NRC is just
a basis of a $q$-dimensional projective subspace.

If $\#F\geq n+2$ or $n=2$, then each automorphic collineation of the
NRC (\ref{param}) preserves osculating subspaces. Otherwise, there
are automorphic collineations of the NRC that do not preserve all
osculating subspaces, whence the concept of osculating subspaces
depends on the parametric representation of the NRC rather than on
the points of the NRC \cite{havl-84}, \cite[2.4]{havl-85}.
\end{rem}

\subsection{Number and dimensions of nuclei}

The following theorem links nuclei of a NRC with Pascal's triangle:

\begin{thm} \label{N-basis}{\rm \cite{gmai+h-00a}}
   If $\# F \geq k+1$, then the $k$-nucleus $\ncal{k}\Vcal_1^n$ of
   the normal rational curve (\ref{param}) equals the subspace
   spanned by those base points $F \cbf_j$,
   where $j\in\{0,1,\ldots,n\}$ is subject to
   \begin{equation}\label{N-basis-j}
   {k+1 \choose j}\equiv {k+2 \choose j}\equiv \ldots\equiv {n \choose j}
   \equiv 0 \mod{\chara F}.
   \end{equation}
\end{thm}
By Theorem \ref{N-basis}, $\chara F=0$ implies that all nuclei of a
NRC are empty. Thus we assume in the remaining part of this section
that
   \begin{equation}
   \chara F =:  p > 0;\;
   n =:  \langle n_\sigma\rangle,\;
   n+1      =:  b=:\langle b_\sigma\rangle\mbox{ (in base }p).
   \end{equation}

We are now in a position to describe the (projective) dimension of a
$k$-nucleus:

\begin{thm}\label{hauptsatz}{\rm \cite{gmai+h-00a}}
   If $\# F\geq k+1$ and
   \begin{equation}\label{hauptsatz-k}
   T(R,b)=\sum_{\mu=R}^{\infty} b_\mu p^\mu \leq k+1 <
   \sum_{\sigma=Q}^{\infty} b_\sigma p^\sigma = T(Q,b)
   \end{equation}
   with at most one $b_\sigma\neq 0$ for $\sigma\in\{Q,Q+1,\ldots,R-1\}$,
   then the $k$-nucleus of $\Vcal_1^n$ has dimension
   \begin{equation}\label{hauptsatz-dim}
   n - \Big(1 + \sum_{\mu=0}^{R-1} n_\mu p^\mu\Big)
   \prod_{\sigma=R}^\infty (n_\sigma +1) = \Sigma(R,n) - 1.
   \end{equation}
\end{thm}
 The condition on the digits $b_\sigma$ guarantees that the top
line function $T$ does not assume a value that is properly between
$T(R,b)$ and $T(Q,b)$.

If $k=n-1$, then (\ref{hauptsatz-dim}) turns into
Timmermann's formula \cite[4.15]{timm-78}
   \begin{equation}\label{hauptsatz-Timm}
   \dim\ncal{n-1}\Vcal_1^n = n-\prod_{\sigma=0}^\infty (n_\sigma +1)
    = \Sigma(1,n)-1;
   \end{equation}
cf.\ also \cite{timm-77}.

\begin{rem} If the ground field $F$ does not meet the richness
condition of Theorem \ref{hauptsatz}, then (\ref{hauptsatz-dim}) is a
lower bound for the dimension of the $k$-nucleus, but it seems to be
an open problem to explicitly determine the dimension of the
$k$-nucleus in terms of $k$, $n$, and $\#F$. See also \cite{havl-84a}
and Example \ref{beispiel-vero}.
\end{rem}
Next we state a formula for the number of distinct nuclei:

\begin{thm}\label{anzahlsatz}{\rm \cite{gmai+h-00a}}
 If $\#F\geq n$, then there are as many distinct nuclei
 of\/ $\Vcal_1^n$ as non-zero digits in the representation of $b=n+1$
 in base $p$.
\end{thm}

\begin{exa}
Let $p = 2$ and $n = 50$: The representation of $50+1$ in base $2$ is
$\langle 110011\rangle$; there are four non-zero digits. So there are
four distinct nuclei, including one empty nucleus. From
\begin{eqnarray}\label{top-bsp}
   & T(0,51) = \langle 110011\rangle = 51, \quad
     T(1,51) = \langle 110010\rangle = 50, & \nonumber \\
   & T(2,51)=T(3,51)=T(4,51)  =  \langle 110000\rangle
     = 48, & \nonumber\\
   & T(5,51) = \langle 100000\rangle = 32,\quad
     T(6,51) = \langle 0\rangle = 0, &\nonumber
\end{eqnarray}
we obtain the following values:
\begin{center}
\begin{tabular}{|c|c|c|c|c|c|}
  \hline\rule{0mm}{1em}
  $k$  & $-1,0,\ldots,30$ & $31,32,\ldots,46$ & $47,48$ & $49$\\
  \hline\rule{0mm}{1em}
  $\dim\ncal{k}\Vcal_1^n$
       & $-1$             & $12$              & $38$    & $42$\\
  \hline
\end{tabular}
\end{center}
\end{exa}
\begin{rem}
Let $\#F\geq n$. Then
      \begin{equation}\label{N-punkt}
      n = 2p^i - 2\geq 2\mbox{ (for some }i\in\N)
      \end{equation}
is necessary and sufficient for the smallest non-empty nucleus to be
a single point. In fact, this point is $F\cbf_{p^i-1}$.  In
particular, if $F$ is a finite field of even order $q$ then this
point together with $\Vcal_1^n$ is a $(q+2)$-arc provided that $n=2$
or provided that $n=q-2\geq 2$ \cite{thas-69}. Cf.\ also
\cite{glyn-86}, \cite{stor+t-94}. In general, however, the geometric
meaning of this point seems to be unknown.

From Theorem \ref{anzahlsatz} all nuclei are empty exactly if
      \begin{equation}\label{N-leer}
      n = \langle n_J,p-1,\ldots,p-1\rangle = n_Jp^J - 1\geq 2
      \end{equation}
with $1\leq n_J<p$ and $J\in\N$. See also \cite{herz-82},
\cite{karz-87}.

Further properties of nuclei can be found in \cite{gmai-99},
\cite{gmai+h-00a}.
\end{rem}


\subsection{Invariant subspaces}
Each NRC $\Vcal_1^n$ admits a group of collineations that is similar
- via (\ref{g-koo-NRC}) -  to  P$\Gamma$L$(2,F)$ acting on the
projective line $\Pcal(\Xbf)$. If $\#F\geq n+2$ or $n=2$ then this
group is the full collineation group of the curve \cite{havl-84}.

Each nucleus is an {\em invariant subspace} i.e.\ it remains fixed
(as a point set) under the full collineation group of the NRC. In
many low-dimensional examples there are no invariant subspaces other
than nuclei. Clearly, all invariant subspaces form a lattice with the
operations of ``join'' and ``meet''.

In order to find all invariant subspaces, we follow J.\ Gmainer
\cite{gmai-01a}: Suppose that the dimension $n$ is fixed. For $j \in
\N$ let
     \begin{equation} \label{omega-eq}
       \Omega(j):= \{
             m \in \N \mid 0 \leq m \leq n,
             {\textstyle{m \choose j}} \not\equiv 0 \mod{\chara F}
             \}.
     \end{equation}
Given a subset $J \subset \{0,1, \dots, n\}$ then put
\begin{equation}\label{omega-psi}
 \Omega(J):= \bigcup\limits_{j \in J} \Omega(j),\quad
  \Psi (J):= \bigcup\limits_{j \in J} \{j, n-j \}.
\end{equation}
Both $\Omega$ and $\Psi$ are closure operators on $\{0,1,\ldots,n\}$.
The idea behind these definitions is as follows: For each $a\in
F\setminus \{0\}$ the mapping
\begin{equation}\label{diagonal}
  F(x_0,x_1,\ldots,x_n)\mapsto F(a^0 x_0,a^1 x_1,\ldots,a^n x_n)
\end{equation}
describes (in terms of coordinates) an automorphic projective
collineation of the NRC (\ref{param}). If the ground field is
sufficiently large then every subspace, which is invariant under all
collineations (\ref{diagonal}), is spanned by base points
$F\cbf_\lambda$. Here $\lambda$ runs in a subset $\Lambda$ of
$\{0,1,\ldots,n\}$. So it is enough to look for appropriate subsets
$\Lambda$. Likewise, the projective collineation
\begin{equation}
  F(x_0,x_1,\ldots,x_n)\mapsto F(x_n,x_{n-1},\ldots,x_0)
\end{equation}
leaves the NRC invariant, whence $\Lambda$ has to be closed with
respect to $\Psi$. A similar argument yields the operator $\Omega$.

Now we are able to formulate the main theorem for invariant
subspaces.

\begin{thm}  \label{main-theo}{\rm \cite{gmai-01a}}
 Let $\#F\geq n+2$ or $n=2$.
 A subspace $\Ucal$ is invariant under the collineation group of the
 normal rational curve (\ref{param}) if, and only if, $\Ucal$ is
 spanned by base points $F\cbf_\lambda$ with $\lambda \in
 \Lambda\subset \{0,1, \dots, n \}$ such that $\Psi (\Lambda)
 \subset\Lambda$ and $\Omega(\Lambda) \subset \Lambda$.
\end{thm}
In case of $\chara F= 0$ there are only trivial invariant
subspaces. Thus we may restrict ourselves to the case
\begin{equation}
  \chara F = p>0.
\end{equation}
By Theorem \ref{main-theo} it suffices to find all $\Psi$-closed
index sets $\Omega(J)\subset\{0,1,\ldots,n\}$. To this end we proceed
in four steps:

Firstly, let $\langle b_\sigma\rangle$ be the expansion of $b=n+1$ in
base $p$. We define
\begin{equation}
            V(i,b):=\sum_{\sigma=0}^{i-1} b_\sigma p^\sigma
            \mbox{ for all }i\in\N.
\end{equation}

Secondly, we fix one number $i\in\N$. Suppose that $I_\alpha$, where
$\alpha\in \{1,2,\ldots,L\}$, is a family of sets such that the
following conditions on the sets $I_\alpha$ and the digits $b_\sigma$
of $b$ hold true:
 \begin{enumerate}\itemsep0em\parsep0em
  \item Each non-empty set $I_\alpha$ has the form
   $I_\alpha=\{j\in\N\mid H_\alpha > j\geq h_\alpha\}$ with
   $i-1\geq H_\alpha>h_\alpha\geq 0$.
  \item If $\alpha>\beta$ and $I_\alpha,I_\beta\neq\emptyset$,
   then $h_\alpha>H_\beta$.
  \item $b_{H_\alpha}<p-1$ and $b_{h_\alpha}>0$ for each non-empty
  set $I_\alpha$.
  \end{enumerate}
For empty subsets $I_\alpha$ no numbers $H_\alpha,h_\alpha$ will be
defined. Thus
\begin{equation}\label{V(i,b)}
    V(i,b)= \langle \ldots,
    b_{H_\alpha},\underbrace{b_{H_\alpha-1},\ldots,b_{h_\alpha}}_{I_\alpha\neq\emptyset},
    \ldots\rangle
\end{equation}
and blocks of digits belonging to different non-empty sets
$I_\alpha\cup\{H_\alpha\}$ do not overlap. So we are in a position to
define a number by simultaneously changing the digits of $V(i,b)$ for
all non-empty sets $I_\alpha$ as follows:
\begin{equation}\label{V(I,i,b)}
\begin{array}{l}
  V(I_1, \dots, I_L; i,b):=
  \langle \ldots,
   b_{H_\alpha}+1,\underbrace{0, \ldots,0}_{I_\alpha\neq\emptyset},
  \dots\rangle
\end{array}
\end{equation}

Thirdly, we assign to each $(I_1,I_2,\ldots,I_L,i,b)$, such that
$V(I_1,I_2,\ldots,I_L;i,b)$ is defined, the set
\begin{equation}\label{Tcal}
   \Tcal(I_1\times I_2\times\cdots \times I_L)
\end{equation}
of all $(T_1,T_2,\ldots,T_L)$ satisfying the following conditions:
 \begin{enumerate}\itemsep0em\parsep0em
  \item If $T_\alpha\neq\emptyset$ then $T_\alpha\subset I_\alpha$
        and $h_\alpha=\min T_\alpha$.
  \item $V(T_1,T_2,\ldots,T_L;i,b)$ is defined.
  \end{enumerate}

Finally, whenever $V(I_1,I_2,\ldots,I_L;i,b)$ is defined, put
  \begin{equation} \label{union-sets}
      \Lambda (I_1, \dots, I_L;i,b):=
                 \bigcup \Big( \Omega ( V(T_1, \dots, T_L; i,b) )
                         \Big),
  \end{equation}
by taking the union over all $(T_1,T_2,\ldots,T_L) \in
\Tcal(I_1\times I_2\times\cdots\times I_L)$.

Let us say that an invariant subspace is {\em irreducible} if it is
not spanned by the invariant subspaces properly contained in it.
Then, with all the assumptions made so far, we obtain the following
result:

\begin{thm} \label{all irreducible}{\rm \cite{gmai-01a}}
 Let $\#F\geq n+2$ or  $n=2$. An
invariant subspace $\Ucal$ of the normal rational curve (\ref{param})
is irreducible if, and only if, it can be written as
\begin{equation}
   \Ucal:= {\mathrm {span}} \{F\cbf_\lambda \mid \lambda \in
                 \Lambda (I_1, \dots, I_L;i,b) \}.
\end{equation}
\end{thm}
 As the lattice of invariant subspaces has only finitely many
elements, each invariant subspace is a join of irreducible ones.

\begin{exa}
Let $n=31$, $p=3$ and $\# F\geq 33$. From $b=32=\langle 1012\rangle$
we get
\begin{displaymath}
\begin{array}{rcl@{\hspace{1em}}rcl}
  V(0,32)         &=& \langle 0 \rangle,&
  V(1,32)         &=& \langle 2 \rangle, \\
  V(3,32)         &=& \langle 012 \rangle,&
  V(\{0\};3,32)   &=& \langle 020 \rangle,\\
  V(\{1\};3,32)   &=& \langle 102 \rangle,&
  V(\{0,1\};3,32) &=& \langle 100 \rangle,\\
  V(4,32) &=&\langle 1012 \rangle.&&
\end{array}
\end{displaymath}
Note that, for example, $V(0,32)=V(\emptyset;0,32)$. Further $V(2,32)
= V(3,32)$, $V(\{0\};2,32) = V(\{0\};3,32)$, and
$V(I_1,I_2,\ldots,I_L;4,32) \geq 32$. So
\begin{displaymath}
\begin{array}{r@{{}={}}l}
 \Omega(\langle 0 \rangle) & \{\langle 0 \rangle,
                                 \langle 1 \rangle,\ldots,
                                 \langle 1011\rangle\}, \\
 \Omega(\langle 2 \rangle) & \{\langle 2 \rangle,
                                 \langle 12 \rangle,
                                 \langle 22 \rangle,
                                 \langle 102 \rangle,
                                 \langle 112 \rangle,
                                 \langle 122 \rangle,
                                 \langle 202 \rangle,
                                 \langle 212 \rangle,
                                 \langle 222 \rangle,
                                 \langle 1002 \rangle
                                 \}, \\
 \Omega(\langle 12\rangle) & \{\langle 12 \rangle,
                                 \langle 22 \rangle,
                                 \langle 112 \rangle,
                                 \langle 122 \rangle,
                                 \langle 212 \rangle,
                                 \langle 222 \rangle
                                 \}, \\
 \Omega(\langle 20\rangle) &  \{\langle 20 \rangle,
                                  \langle 21 \rangle,
                                  \langle 22 \rangle,
                                  \langle 120 \rangle,
                                  \langle 121 \rangle,
                                  \langle 122 \rangle,
                                  \langle 220 \rangle,
                                  \langle 221 \rangle,
                                  \langle 222 \rangle
                                 \}, \\
 \Omega(\langle 102\rangle)& \{\langle 102 \rangle,
                                 \langle 112 \rangle,
                                 \langle 122 \rangle,
                                 \langle 202 \rangle,
                                 \langle 212 \rangle,
                                 \langle 222 \rangle
                                 \}, \\
 \Omega(\langle 100\rangle)& \{\langle 100 \rangle,
                                 \langle 101 \rangle,\ldots,
                                 \langle 222 \rangle
                                 \}, \\
 \Omega(\langle 1012\rangle)& \emptyset,
\end{array}
\end{displaymath}
are the relevant index sets and
\begin{displaymath}
\begin{array}{r@{\:{}={}\:}l}
 \Lambda(\emptyset;0,32) &\Omega(\langle 0\rangle), \\
 \Lambda(\emptyset;1,32) & \Omega(\langle 2\rangle), \\
 \Lambda(\emptyset;3,32) & \Omega(\langle 12\rangle), \\
 \Lambda(\{0\};3,32)     & \Omega(\langle 12\rangle)\cup
                             \Omega(\langle 20\rangle)\\
 \Lambda(\{1\};3,32)     & \Omega(\langle 12\rangle)\cup
                             \Omega(\langle 102\rangle)\\
 \Lambda(\{0,1\};3,32)   & \Omega(\langle 12\rangle)\cup
                             \Omega(\langle 20\rangle)\cup
                             \Omega(\langle 102\rangle)\cup
                             \Omega(\langle 100\rangle)\\
 \Lambda(\emptyset;4,32) & \emptyset.
\end{array}
\end{displaymath}
The Hasse diagram of the lattice of invariant subspaces is given in
the figure. Filled circles represent irreducible subspaces and double
circles mark nuclei.

  {\unitlength0.7cm
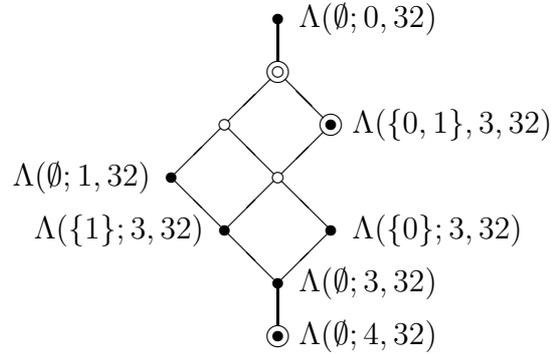
    \begin{center}
      \begin{picture}(6.0,6.0)
 \put(2,6){\circle*{0.2}}
 \put(2.4,6){\makebox(0,0)[ml]{$\Lambda(\emptyset;0,32)$}}

 \put(2,5.2){\line(0,1){0.7}}

 \put(2,5){\circle{0.2}}
 \put(2,5){\circle{0.4}}
 \put(2.4,5){\makebox(0,0)[ml]{}}

 \put(1.07,4.07){\line( 1,1){0.79}}
 \put(2.86,4.14){\line(-1,1){0.72}}

 \put(1,4){\circle{0.2}}
 \put(3,4){\circle*{0.2}}\put(3,4){\circle{0.4}}
 \put(3.4,4){\makebox(0,0)[ml]{$\Lambda(\{0,1\},3,32)$}}

 \put(0.07,3.07){\line( 1,1){0.86}}
 \put(1.93,3.07){\line(-1,1){0.86}}
 \put(2.07,3.07){\line( 1,1){0.79}}

 \put(0,3){\circle*{0.2}}
 \put(-0.4,3){\makebox(0,0)[mr]{$\Lambda(\emptyset;1,32)$}}
 \put(2,3){\circle{0.2}}

 \put(0.93,2.07){\line(-1,1){0.86}}
 \put(1.07,2.07){\line( 1,1){0.86}}
 \put(2.93,2.07){\line(-1,1){0.86}}

 \put(1,2){\circle*{0.2}}
 \put(0.6,2){\makebox(0,0)[mr]{$\Lambda(\{1\};3,32)$}}
 \put(3,2){\circle*{0.2}}
 \put(3.4,2){\makebox(0,0)[ml]{$\Lambda(\{0\};3,32)$}}

 \put(1.93,1.07){\line(-1,1){0.86}}
 \put(2.07,1.07){\line(1,1){0.86}}

 \put(2,1){\circle*{0.2}}
 \put(2.4,1){\makebox(0,0)[ml]{$\Lambda(\emptyset;3,32)$}}

 \put(2,0.2){\line(0,1){0.7}}

 \put(2,0){\circle*{0.2}}\put(2,0){\circle{0.4}}
 \put(2.4,0){\makebox(0,0)[ml]{$\Lambda(\emptyset;4,32)$}}
       \end{picture}
    \end{center}
}
\end{exa}
In many low-dimensional examples the invariant subspaces form a
chain. In general, however, the following holds:

\begin{thm}{\rm \cite{gmai-01a}}
      Let the positions of the non-zero digits of $b:=n+1$ in base $p$ be denoted
      by $N_1, N_2, \dots, N_d$. Then the lattice of
      invariant subspaces is totally ordered if, and only
      if, one of the following cases occurs:
        \begin{enumerate}\itemsep0em\parsep0em
         \item $d \in \{1,2 \}$.
         \item $d \ge 3,\, N_d - N_1 = d-1,$ and $N_2 = \dots = N_{d-1} =
               p-1$.
        \end{enumerate}
\end{thm}
Thus all invariant subspaces can be found, provided that the ground
field is sufficiently large. Also, in specific cases the structure of
the lattice of invariant subspaces is known.

\begin{rem}
If we project a NRC from one of its invariant subspaces other than
$\Pcal(\Ybf)$, then a rational curve is obtained; this curve admits a
collineation group isomorphic to P$\Gamma$L$(2,F)$. Via
(\ref{g-koo-NRC}) and the projection, the group actions on the curve
and the projective line $\Pcal(\Xbf)$ are similar.
\end{rem}

\section{Pascal's simplex modulo a prime}\label{Kap-Multinom}

Throughout this section let $p$ be a fixed prime. Given
$m,t\in\N$ then put
   \begin{equation}\label{E-m-t}
   {E^t_m} :=\{ (e_0,e_1,\ldots,e_m)\in\N^{m+1} \mid e_0+e_1+\ldots+e_m=t
   \}.
   \end{equation}
The array of multinomial coefficients ${t \choose
{e_0,e_1,\ldots,e_m}}$ with $(e_0,e_1,\ldots,e_m)\in{E^t_m}$
is frequently called {\em Pascal's simplex}.

The theorem of Lucas (\ref{lucas}) can be generalized to multinomial
coefficients as follows \cite[364]{brou+w-95}: If
$t,e_0,e_1,\ldots,e_m\in\N$ have re\-pre\-sent\-ations $t=\sum_\sigma
t_{\sigma}p^\sigma$ and $e_i=\sum_\sigma e_{i,\sigma}p^\sigma$ in
base $p$ then
   \begin{equation}\label{lucas-multi}
   {t \choose {e_0,e_1,\ldots,e_m}}\equiv
   \prod_{\sigma\in\N} {t_\sigma \choose
   {e_{0,\sigma},e_{1,\sigma},\ldots,e_{m,\sigma}}}
   \mod p.
   \end{equation}
For trinomial coefficients ($m=3$) it is possible to illustrate
Pascal's simplex in the form of a pyramid:
    
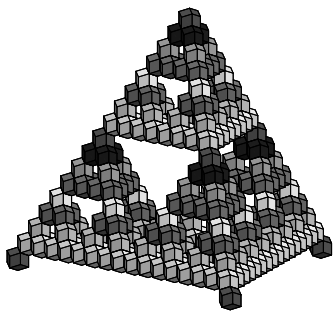
\begin{center}
         {\includegraphics[height=4.0cm]{py2-16.eps}}
    \end{center}
The picture above shows a part of Pascal's pyramid modulo $2$. It is
based upon a tiling of the space by rhombic dodecahedra. If an entry
of the pyramid vanishes, then the corresponding dodecahedron is
omitted. Entries at the same ``horizontal'' level ($t$ constant) are
equally shaded. Cf.\ \cite{hilt+p-99}.

The following has been established independently by F.T.\ Howard
\cite[Theorem 3.1]{howa-74} and N.A.\ Volodin \cite[Theorem
2]{volo-89}; see also \cite{volo-94}:

The number of $(m+1)$--tuples $(e_0,e_1,\ldots,e_m)\in {E^t_m}$ such
that the multinomial coefficient ${t \choose {e_0,e_1,\ldots,e_m}}$
is divisible by the prime $p$ equals
   \begin{equation}\label{anzahl}
   {m+t \choose t} - \prod_{\sigma\in\N}{{m+t_\sigma}\choose t_\sigma}.
   \end{equation}

\section{Veronese varieties}\label{Kap-Vero}
\subsection{Definition of $(r,k)$-nuclei}

Let $\{\bbf_0,\bbf_1,\ldots,\bbf_m\}$ be a basis of an
$(m+1)$-dimensional vector space $\Xbf$ over $F$ (the parameter
space) and let $\Ybf$ be an ${m+t}\choose t$-dimensional vector space
over $F$ with a basis $\{\cbf_{e_0,e_1,\ldots,e_m}\mid
(e_0,e_1,\ldots,e_m)\in {E^t_m}\}$; cf.\ (\ref{E-m-t}). We shall
always assume that $m\geq 1$ and $t\geq 2$ in order to avoid
trivialities.

Generalizing (\ref{g-koo-NRC}), the {\em Veronese mapping} is given
by
   \begin{equation}\label{g-koo}
   F \Big( \sum_{i=0}^m x_i \bbf_i \Big)
   \mapsto
   F \Big( \sum_{{E^t_m}} x_0^{e_0} x_1^{e_1} \ldots x_m^{e_m}
   \cbf_{e_0,e_1,\ldots,e_m} \Big)
   \;\;(x_i\in F).
   \end{equation}
Its image is a {\em Veronese variety} $\Vcal_m^t$ with ambient space
$\Pcal(\Ybf)$, i.e.\ the projective space on $\Ybf$. (By putting
$m:=1$ and $n:=t$ a NRC $\Vcal_1^n$ is obtained.)

The Veronese image of each $r$-dimensional subspace of $\Pcal(\Xbf)$
$(0\leq r< m)$ is a sub-Veronesean $\Vcal_r^t$ of $\Vcal_m^t$. (For
$r=0$ we get just one point, for $r=1$ a normal rational curve, etc.
Cf.\ also \cite{bura-74}.) For each $k\in\{-1,0,\ldots,t-1\}$ there
exists a $k$-{\em osculating subspace of $\Vcal_m^t$ along}
$\Vcal_r^t$. We call it an {\em $(r,k)$-osculating subspace} of
$\Vcal_m^t$. Its dimension equals
\begin{equation}\label{schmieg-dim}
    \sum_{i=t-k}^{t}{r+i \choose i}{m+t-r-i-1 \choose t-i} -1;
\end{equation}
cf.\ \cite{herz-82} and the papers cited in Remark \ref{def-osk}. We
are thus led to the following definition:
\begin{defi}
The $(r,k)$-{\em nucleus} of a Veronese variety $\Vcal_m^t$ is the
intersection of all its $(r,k)$-osculating subspaces.
\end{defi}
The $k$-nuclei of a normal rational curve are the $(0,k)$-nuclei
according to the present definition.
\begin{rem}
A geometric characterization of quadratic Veronese mappings ($t=2$)
can be found in \cite{havl+z-97}. Combinatorial characterizations of
the Veronese surface ($m=t=2$) over a finite field and further
references on this particular subject are given in \cite{hirs+t-91}.
Applications of Veronese varieties over finite fields in coding
theory and authentication systems can be found in \cite{gopp-88},
\cite{havl-98}, \cite{havl-99}, \cite{hirs+s-98}, \cite{stor+t-94},
\cite{thas-92}, \cite{zane-98}. Partial linear spaces derived from
Veronese varieties are discussed in \cite{melo-83}.
\end{rem}

\subsection{Intersection of osculating hyperplanes}

From (\ref{schmieg-dim}), each $(t-1,m-1)$-osculating subspace of a
Veronese variety $\Vcal_m^t$ is a hyperplane of $\Pcal(\Ybf)$ which
is called an {\em osculating hyperplane} (or {\em contact
hyperplane}) of the Veronese variety $\Vcal_m^t$. Thus to each
hyperplane of the parameter space there corresponds an osculating
hyperplane of the Veronesean. In terms of dual bases this {\em dual
Veronese mapping} is given by
   \begin{equation}\label{symm-potenz}
   F\Big(\sum_{i=0}^m a_i \bbf_i^\ast\Big)
   \mapsto
   F \Big(
   \sum_{{E^t_m}} {\textstyle {t \choose {e_0,e_1,\ldots,e_m}}}
   a_0^{e_0} a_1^{e_1}\ldots a_m^{e_m}
   \cbf_{ {e_0},{e_1}\ldots,{e_m} }^\ast
    \Big)
   \;\;(a_i\in F).
   \end{equation}
See also \cite[pp.\ 160--163]{bura-61}. The intersection of all
osculating hyperplanes of a $\Vcal_m^t$ is its $(m-1,t-1)$-nucleus.
Both A.\ Herzer \cite{herz-82} and H.\ Karzel \cite{karz-87}
determined all Veronese varieties where this specific nucleus is
empty.

\begin{thm}\label{theo-1}{\rm \cite{gmai+h-00b}}
 The $(m-1,t-1)$-nucleus of a Veronese
 variety $\Vcal_m^t$ contains exactly those base points
 $F\cbf_{e_0,e_1,\ldots,e_m}$ satisfying
   \begin{equation}\label{verschwindet}
   {t \choose {e_0,e_1,\ldots,e_m}} \equiv 0 \mod {\chara F}.
   \end{equation}
 If $\# F\geq t$, then this nucleus is spanned by those base
 points.
 \end{thm}
From this and (\ref{anzahl}) follows

\begin{thm}\label{theo-2}{\rm \cite{gmai+h-00b}}
 Let $\sum_{\sigma\in\N} t_\sigma p^\sigma$ be the representation
 of $t$ in base $p=\chara F > 0$. If $\# F\geq t$,
 then the $(m-1,t-1)$-nucleus of a Veronese variety $\Vcal_m^t$ has
 dimension
   \begin{equation}\label{knotendim}
   {{m+t}\choose t} -
   \prod_{\sigma\in\N} {{m+t_\sigma}\choose t_\sigma}-1.
   \end{equation}
\end{thm}
\begin{exa}\label{beispiel-vero}
Let $\chara F=2$.

The $(1,1)$-nucleus of the Veronese surface $\Vcal_2^2$ is a plane;
cf.\ \cite[Chapter 25]{hirs+t-91}.

From Theorem \ref{theo-2} the $(1,2)$-nucleus of the Veronese surface
$\Vcal_2^3$ is a single point provided that $\#F\neq 2$. On the other
hand, if $\#F=2$ then, by solving a system of seven linear equations,
the $(1,2)$-nucleus of $\Vcal_2^3$ is easily seen to be
three-dimensional. In either case $\Vcal_2^3$ carries a family of
twisted cubics that arise as Veronese images of the lines in the
parameter plane. For $\#F\neq 2$ the $2$-nucleus of a twisted cubic
is empty, but for $\#F=2$ this nucleus is a single point.
\end{exa}

{\small

\begin{thebibliography}{10}\itemsep-3pt

\bibitem{bert-07}
E.~Bertini.
\newblock {\em Introduzione alla geometria proiettiva degli iperspazi}.
\newblock E.\ Spoerri, Pisa, 1907.

\bibitem{bert-24}
E.~Bertini.
\newblock {\em Einf\"uhrung in die projektive {G}eometrie mehrdimensionaler
  {R}\"au\-me}.
\newblock Seidel u.\ Sohn, Wien, 1924.

\bibitem{brau-762}
H.~Brauner.
\newblock {\em Geometrie projektiver R\"aume II}.
\newblock BI-Wissenschaftsverlag, Mannheim Wien Z\"urich, 1976.

\bibitem{brou+w-95}
A.E.\ Brouwer and H.A.\ Wilbrink.
\newblock Block designs.
\newblock In F.\ Buekenhout, editor, {\em Handbook of incidence geometry},
  chapter~8, pages 349--382. Elsevier, Amsterdam, 1995.

\bibitem{bura-61}
W.~Burau.
\newblock {\em Mehrdimensionale projektive und h\"ohere {G}eometrie}.
\newblock Dt.\ Verlag d.\ Wissenschaften, Berlin, 1961.

\bibitem{bura-74}
W.~Burau.
\newblock {\"U}ber ausgezeichnete {A}ufspaltungen des {R}aumes einer
  {V}eroneseschen ${V}^t_n$ und ihre {A}nwendung auf die {B}erechnung der
  {H}ilbertfunktion der rationalen {N}ormregelmannigfaltigkeiten.
\newblock {\em Ver\"off.\ Univ.\ Innsbruck}, 91:17--28, 1974.

\bibitem{fine-47}
N.J.\ Fine.
\newblock Binomial coefficients modulo a prime.
\newblock {\em Am.\ Math.\ Mon.}, 54:589--592, 1947.

\bibitem{glyn-86}
D.G.\ Glynn.
\newblock The non-classical $10$-arc of {PG}$(4,9)$.
\newblock {\em Discrete Math.}, 59:43--51, 1986.

\bibitem{gmai-99}
J.~Gmainer.
\newblock {\em Rationale {N}ormkurven in {R}\"aumen mit positiver
  {C}harakteristik}.
\newblock {T}hesis, Vienna University of Technology, 1999.

\bibitem{gmai-01a}
J.~Gmainer.
\newblock Pascal's triangle, normal rational curves, and their invariant
  subspaces.
\newblock {\em Eur.\ J.\ Comb.}, 22:37--49, 2001.

\bibitem{gmai+h-00b}
J.\ Gmainer and H.~Havlicek.
\newblock A dimension formula for the nucleus of a {V}eronese variety.
\newblock {\em Lin.\ Algebra Appl.}, 305:191--201, 2000.

\bibitem{gmai+h-00a}
J.\ Gmainer and H.~Havlicek.
\newblock Nuclei of normal rational curves.
\newblock {\em J.\ Geometry}, 69:117--130, 2000.

\bibitem{gopp-88}
V.D.\ Goppa.
\newblock {\em Geometry and codes}.
\newblock Kluwer, Dordrecht Boston London, 1988.

\bibitem{harb-75}
H.~Harborth.
\newblock {\"U}ber die {T}eilbarkeit im {P}ascal-{D}reieck.
\newblock {\em Math.-phys.\ Semesterber.}, 22:13--21, 1975.

\bibitem{hass-37}
H.~Hasse.
\newblock Noch eine {B}egr\"undung der {T}heorie der h\"oheren
  {D}ifferentialquotienten in einem algebraischen {F}unktionenk\"orper einer
  {U}nbestimmten.
\newblock {\em J.\ Reine Angew.\ Math.}, 177:215--237, 1937.

\bibitem{havl-83}
H.~Havlicek.
\newblock Normisomorphismen und {N}ormkurven endlichdimensionaler projektiver
  {D}esargues-{R}\"aume.
\newblock {\em Monatsh.\ Math.}, 95:203--218, 1983.

\bibitem{havl-84}
H.~Havlicek.
\newblock Die automorphen {K}ollineationen nicht entarteter {N}ormkurven.
\newblock {\em Geom.\ Dedicata}, 16:85--91, 1984.

\bibitem{havl-84a}
H.~Havlicek.
\newblock Erzeugnisse projektiver {B}\"undelisomorphismen.
\newblock {\em Ber.\ Math.-Stat.\ Sekt.\ For\-schungs\-zent.\ Graz}, 215, 1984.

\bibitem{havl-85}
H.~Havlicek.
\newblock Applications of results on generalized polynomial identities in
  desarguesian projective spaces.
\newblock In R.\ Kaya, P.\ Plaumann, and K.\ Strambach, editors, {\em Rings and
  geometry}, pages 39--77. D.\ Reidel, Dordrecht, 1985.

\bibitem{havl-98}
H.~Havlicek.
\newblock The {V}eronese surface in {PG}$(5,3)$ and {W}itt's
  $5$-$(12,6,1)$--design.
\newblock {\em J.\ Comb.\ Theory, Ser.\ A}, 84(1):87--94, 1998.

\bibitem{havl-99}
H.~Havlicek.
\newblock {G}iuseppe {V}eronese and {E}rnst {W}itt -- neighbours in {PG}(5,3).
\newblock {\em Aequationes Math.}, 58:85--92, 1999.

\bibitem{havl+z-97}
H.\ Havlicek and C.\ Zanella.
\newblock Quadratic embeddings.
\newblock {\em Beitr.\ Algebra Geom.}, 38:289--298, 1997.

\bibitem{herz-82}
A.~Herzer.
\newblock Die {S}chmieghyperebenen an die {V}eronese-{M}annigfaltigkeit bei
  beliebiger {C}harakteristik.
\newblock {\em J.\ Geom.}, 18:140--154, 1982.

\bibitem{hexe+s-78}
E.\ Hexel and H.~Sachs.
\newblock Counting residues modulo a prime in {P}ascal's triangle.
\newblock {\em Indian J.\ Math.}, 20:91--105, 1978.

\bibitem{hilt+p-99}
P.\ Hilton and J.~Pedersen.
\newblock Relating geometry and algebra in the pascal triangle, hexagon,
  tetrahedron, and cuboctahedron {P}art 1: {B}inomial coefficients, extended
  binomial coefficients and preparation for further work.
\newblock {\em College Mathematics Journal}, 30(3):170--186, 1999.

\bibitem{hirs-85}
J.W.P.\ Hirschfeld.
\newblock {\em Finite projective spaces of three dimensions}.
\newblock Oxford University Press, Oxford, 1985.

\bibitem{hirs-98}
J.W.P.\ Hirschfeld.
\newblock {\em Projective geometries over finite fields}.
\newblock Clarendon Press, Oxford, second edition, 1998.

\bibitem{hirs+s-98}
J.W.P.\ Hirschfeld and L.~Storme.
\newblock The packing problem in statistics, coding theory, and finite
  projective spaces.
\newblock {\em J.\ Stat.\ Plann.\ Inference}, 72(1--2):355--380, 1998.

\bibitem{hirs+t-91}
J.W.P.\ Hirschfeld and J.A.\ Thas.
\newblock {\em General Galois geometries}.
\newblock Oxford University Press, Oxford, 1991.

\bibitem{howa-74}
F.T.\ Howard.
\newblock The number of multinomial coefficients divisible by a fixed power of
  a prime.
\newblock {\em Pac.\ J.\ Math.}, 50:99--108, 1974.

\bibitem{kara-96}
V.V.\ Karachik.
\newblock $p$-latin matrices and {P}ascal's triangle modulo a prime.
\newblock {\em Fibonacci Q.}, 34(4):362--372, 1996.

\bibitem{karz-87}
H.~Karzel.
\newblock {\"U}ber einen {F}undamentalsatz der synthetischen algebraischen
  {G}eo\-metrie von {W.~Burau und H.~Timmermann}.
\newblock {\em J.\ Geom.}, 28:86--101, 1987.

\bibitem{long-81b}
C.T.\ Long.
\newblock Pascal's triangle modulo $p$.
\newblock {\em Fibonacci Q.}, 19:458--463, 1981.

\bibitem{melo-83}
N.~Melone.
\newblock Veronese spaces.
\newblock {\em J.\ Geom.}, 20:169--180, 1983.

\bibitem{ries-81}
R.~Riesinger.
\newblock Normkurven in endlichdimensionalen {D}esarguesr\"aumen.
\newblock {\em Geom.\ Dedicata}, 10:427--449, 1981.

\bibitem{robe-57}
J.B.\ Roberts.
\newblock On binomial coefficient residues.
\newblock {\em Canadian J.\ Math.}, 9:363--370, 1957.

\bibitem{stor+t-94}
L.\ Storme and J.A.\ Thas.
\newblock $k$-arcs and dual $k$-arcs.
\newblock {\em Discrete Math.}, 125, No.1--3:357--370, 1994.

\bibitem{thas-69}
J.A.\ Thas.
\newblock Normal rational curves and $(q+2)$--arcs in a {G}alois space
  ${S}_{q-2,q}$ $(q=2^h)$.
\newblock {\em Atti Accad.\ Naz.\ Lincei, VIII.\ Ser., Rend., Cl.\ Sci.\ Fis.\
  Mat.\ Nat.}, 47:249--252, 1969.

\bibitem{thas-92}
J.A.\ Thas.
\newblock M.{D}.{S}.\ codes and arcs in projective spaces: {A} survey.
\newblock {\em Le Matematiche}, 47(2):315--328, 1992.

\bibitem{timm-77}
H.~Timmermann.
\newblock Descrizioni geometriche sintetiche di geometrie proiettive con
  caratteristica $p>0$.
\newblock {\em Ann.\ Mat.\ Pura Appl.\ IV.\ Ser.}, 114:121--139, 1977.

\bibitem{timm-78}
H.~Timmermann.
\newblock {\em Zur {G}eometrie der {V}eronesemannigfaltigkeit bei end\-li\-cher
  {C}ha\-rak\-te\-ri\-stik}.
\newblock {H}abilitationsschrift, Univ.\ Hamburg, 1978.

\bibitem{volo-89}
N.A.\ Volodin.
\newblock Distribution of polynomial coefficients congruent modulo $p^{N}$.
\newblock {\em Math.\ Notes}, 45:195--199, 1989.

\bibitem{volo-94}
N.A.\ Volodin.
\newblock Number of multinomial coefficients not divisible by a prime.
\newblock {\em Fibonacci Q.}, 32(5):402--406, 1994.

\bibitem{wolf-84}
S.~Wolfram.
\newblock Geometry of binomial coefficients.
\newblock {\em Am.\ Math.\ Mon.}, 91:566--571, 1984.

\bibitem{zane-98}
C.~Zanella.
\newblock Linear sections of the finite {V}eronese varieties and
  authentication systems defined using geometry.
\newblock {\em Des.\ Codes Cryptography}, 13(2):199--212, 1998.

\bibitem{zeug-72}
J.~Zeuge.
\newblock Eine geometrische {K}ennzeichnung der {M}annigfaltigkeiten von
  {S}egre und {V}eronese und eine damit zusammenh\"angende ausgezeichnete
  {T}ransformation zwischen gewissen projektiven {R}\"aumen.
\newblock {\em Atti Accad.\ naz.\ Lincei, VIII.\ Ser., Rend., Cl.\ Sci.\ fis.\
  Mat.\ natur.}, 53:531--540, 1972.

\bibitem{zeug-77}
J.~Zeuge.
\newblock Die {S}chmiegr\"aume an die {V}eronesemannigfaltigkeit.
\newblock {\em Mitt.\ Math.\ Ges.\ Hamburg}, 10(5):391--393, 1977.

\end{thebibliography}

}

\noindent Hans Havlicek, Abteilung f\"ur Lineare Geometrie,
Technische Universit\"at, Wiedner Hauptstra{\ss}e 8--10, A-1040 Wien,
Austria.

\noindent Email: {\tt havlicek@geometrie.tuwien.ac.at}

\noindent Web site: {\tt http://www.geometrie.tuwien.ac.at/havlicek}

\end{document}